\documentclass[10pt,a4paper]{article}
\usepackage{amsfonts,amsmath}

\newtheorem{theorem}{Theorem}

\newtheorem{remark}{Remark}

\input{tcilatex}
\begin{document}

\date{}
\title{\textbf{A remark on the existence of entire large and bounded
solutions to a (}$k_{1}$,$k_{2}$\textbf{)-Hessian system with gradient term}}
\author{ \vspace{1mm} \textsf{{Dragos-Patru Covei} } \\
{\small \textit{Department of Applied Mathematics}}\\
{\small \textit{Bucharest University of Economic Studies }}\\
{\small \textit{Piata Romana, 1st district, postal code: 010374, postal
office: 22, Romania}}\\
{\small \textit{e-mail address:coveidragos@yahoo.com}}}
\maketitle

\begin{abstract}
In this paper, we study the existence of positive entire large and bounded
radial positive solutions for the following nonlinear system 
\begin{equation*}
\left\{ 
\begin{array}{l}
S_{k_{1}}\left( \lambda \left( D^{2}u_{1}\right) \right) +a_{1}\left(
\left\vert x\right\vert \right) \left\vert \nabla u_{1}\right\vert
^{k_{1}}=p_{1}\left( \left\vert x\right\vert \right) f_{1}\left(
u_{2}\right) \text{ for }x\in \mathbb{R}^{N}\text{ , } \\ 
S_{k_{2}}\left( \lambda \left( D^{2}u_{2}\right) \right) +a_{2}\left(
\left\vert x\right\vert \right) \left\vert \nabla u_{2}\right\vert
^{k_{2}}=p_{2}\left( \left\vert x\right\vert \right) f_{2}\left(
u_{1}\right) \text{ for }x\in \mathbb{R}^{N}\text{ .}%
\end{array}%
\right. 
\end{equation*}%
Here $S_{k_{i}}\left( \lambda \left( D^{2}u_{i}\right) \right) $ is the $%
k_{i}$-Hessian operator, $a_{1}$, $p_{1}$, $f_{1}$, $a_{2}$, $p_{2}$ and $%
f_{2}$ are continuous functions. Our results give an answer of the question
raised in \cite{GR}.
\end{abstract}

\vspace{0.5cm}{\scriptsize 2010 AMS Subject Classification: Primary: 35J25,
35J47 Secondary: 35J96.}


\section{Introduction}

Explosive or bounded radial solutions of elliptic systems have been studied
intensively in the last few decades (see among others Alves and Holanda \cite%
{CA}, Bandle and Giarrusso \cite{BA}, Cirstea and R\u{a}dulescu \cite{CR}, Cl%
\'{e}ment-Man\'{a}sevich and Mitidieri \cite{CL}, the author \cite{CD2}, De
Figueiredo and Jianfu \cite{FI}, Galaktionov and V\'{a}zquez \cite{GA},\
Ghergu and R\u{a}dulescu \cite{GR}, Lair and Wood \cite{LA}, Lair \cite%
{LAIR2, LAIR3}, Peterson and Wood \cite{PW}, Quittner \cite{Q} and Zhang and
Zhou \cite{ZZ}). Most of these works have studied the follwing system 
\begin{equation}
\left\{ 
\begin{array}{l}
\Delta u_{1}=p_{1}\left( \left\vert x\right\vert \right) f_{1}\left(
u_{2}\right) \text{ for }x\in \mathbb{R}^{N}\text{ (}N\geq 3\text{), } \\ 
\Delta u_{2}=p_{2}\left( \left\vert x\right\vert \right) f_{2}\left(
u_{1}\right) \text{ for }x\in \mathbb{R}^{N}\text{ (}N\geq 3\text{).}%
\end{array}%
\right.  \label{cr}
\end{equation}%
C\^{\i}rstea and R\u{a}dulescu \cite[(Theorem 1, p. 828)]{CR}, proved that
the system (\ref{cr}) has infinitely many positive entire large solutions
provided that $f_{1}$, $f_{2}$ are positive locally H\"{o}lder with exponent 
$\beta \in \left( 0,1\right) $, non-decreasing and satisfying%
\begin{equation}
\lim_{t\rightarrow \infty }\frac{f_{2}\left( af_{1}\left( t\right) \right) }{%
t}=0\text{ for all constants }a>1.  \label{vic}
\end{equation}%
In their work, the weights $p_{1}$, $p_{2}$ are positive, symmetric locally H%
\"{o}lder functions in $\mathbb{R}^{N}$ and such that%
\begin{equation}
\int_{0}^{\infty }r^{1-N}\int_{0}^{r}t^{N-1}p_{1}\left( t\right)
dtdr=\int_{0}^{\infty }r^{1-N}\int_{0}^{r}t^{N-1}p_{2}\left( t\right)
dtdr=\infty .  \label{5}
\end{equation}%
In a subsequent paper, Ghergu and R\u{a}dulescu \cite[(Theorem 3, p. 438)]%
{GR} established the existence of solutions for 
\begin{equation}
\left\{ 
\begin{array}{l}
\Delta u_{1}+\left\vert \nabla u_{1}\right\vert =p_{1}\left( x\right)
f_{1}\left( u_{2}\right) \text{ for }x\in \mathbb{R}^{N}\text{ (}N\geq 3%
\text{), } \\ 
\Delta u_{2}+\left\vert \nabla u_{2}\right\vert =p_{2}\left( x\right)
f_{2}\left( u_{1}\right) \text{ for }x\in \mathbb{R}^{N}\text{ (}N\geq 3%
\text{).}%
\end{array}%
\right.  \label{gr}
\end{equation}%
They replaced the condition (\ref{5}) with%
\begin{equation}
\int_{1}^{\infty }e^{-r}r^{1-N}\int_{0}^{r}e^{t}t^{N-1}p_{1}\left( t\right)
dtdr=\int_{1}^{\infty }e^{-r}r^{1-N}\int_{0}^{r}e^{t}t^{N-1}p_{2}\left(
t\right) dtdr=\infty ,  \label{6}
\end{equation}%
while keeping all other conditions on $p_{1}$, $f_{1}$, $p_{2}$ and $f_{2}$.
Moreover, they noticed although 
\begin{equation*}
f_{1}\left( t\right) =\sqrt{t}\text{, }f_{2}\left( t\right) =t\text{, }%
p_{1}\left( r\right) =4\frac{r^{3}+\left( N+2\right) r^{2}}{\sqrt{r^{2}+1}}%
\text{, }p_{2}\left( r\right) =2\frac{r+N}{r^{4}+1},
\end{equation*}%
doesn't satisfy (\ref{gr}), the corresponding system has the positive entire
large solution $\left( \left\vert x\right\vert ^{4}+1,\left\vert
x\right\vert ^{2}+1\right) $. Therefore it is only natural to weaken the
assumptions (\ref{6}). Inspired by this, our analysis, is developed for the
more general nonlinear systems 
\begin{equation}
\left\{ 
\begin{array}{l}
S_{k_{1}}\left( \lambda \left( D^{2}u_{1}\right) \right) +a_{1}\left(
\left\vert x\right\vert \right) \left\vert \nabla u_{1}\right\vert
^{k_{1}}=p_{1}\left( \left\vert x\right\vert \right) f_{1}\left(
u_{2}\right) \text{ for }x\in \mathbb{R}^{N}\text{ (}N\geq 3\text{), } \\ 
S_{k_{2}}\left( \lambda \left( D^{2}u_{2}\right) \right) +a_{2}\left(
\left\vert x\right\vert \right) \left\vert \nabla u_{2}\right\vert
^{k_{2}}=p_{2}\left( \left\vert x\right\vert \right) f_{2}\left(
u_{1}\right) \text{ for }x\in \mathbb{R}^{N}\text{ (}N\geq 3\text{),}%
\end{array}%
\right.  \label{11}
\end{equation}%
where $a_{1}$, $p_{1}$, $f_{1}$, $a_{2}$, $p_{2}$, $f_{2}$ are continuous
functions satisfying certain monotonicity properties, $k_{1},k_{2}\in
\left\{ 1,2,...,N\right\} $ and $S_{k_{i}}\left( \lambda \left(
D^{2}u_{i}\right) \right) $ stands for the $k_{i}$-Hessian operator defined
as the sum of all $k_{i}\times k_{i}$ principal minors of the Hessian matrix 
$D^{2}u_{i}$. The following well known operators 
\begin{equation*}
\begin{tabular}{llll}
\hline
\textbf{Operator:} & Laplacian &  & Monge--Amp\`{e}re \\ \hline
& $S_{1}\left( \lambda \left( D^{2}u_{i}\right) \right) =\Delta u_{i}=\func{%
div}\left( \nabla u_{i}\right) $ &  & $S_{N}\left( \lambda \left(
D^{2}u_{i}\right) \right) =\det \left( D^{2}u_{i}\right) ,$ \\ \hline
\end{tabular}%
\end{equation*}%
are special $k_{i}$-Hessian operator.

The system (\ref{11}) has been the subject of rather deep investigations
since it appears in many branches of applied mathematics (for more on this
see the papers of Bao-Ji and Li \cite{BAOII}, Salani \cite{SA}, Ji and Bao 
\cite{BAO}, Viaclovsky \cite{V,VI} and Zhang and Zhou \cite{ZZ}).

\section{The main results}

We work under the following assumptions:

(P1)\quad $a_{1},a_{2}:\left[ 0,\infty \right) \rightarrow \left[ 0,\infty
\right) $ and $p_{1},p_{2}:\left[ 0,\infty \right) \rightarrow \left(
0,\infty \right) $ are spherically symmetric continuous functions (i.e.,%
\textit{\ }$p_{i}\left( x\right) =p_{i}\left( \left\vert x\right\vert
\right) $ and $a_{i}\left( x\right) =a_{i}\left( \left\vert x\right\vert
\right) $ for $i=1,2$)\textit{;}

(C1)\quad $f_{1}$, $f_{2}:\left[ 0,\infty \right) \rightarrow \left[
0,\infty \right) $ are continuous, increasing, $f_{i}\left( 0\right) \geq 0$
and $f_{i}\left( s\right) >0$ for all $s>0$ with $i=1,2$;

(C2)\quad there exist positive constants $\overline{c}_{1},\overline{c}_{2}$%
, the continuous and increasing functions $h_{1}$, $h_{2}:\left[ 0,\infty
\right) \rightarrow \left[ 0,\infty \right) $ and the continuous functions $%
\overline{\varphi }_{1}$, $\overline{\varphi }_{2}:\left[ 0,\infty \right)
\rightarrow \left[ 0,\infty \right) $ such that%
\begin{eqnarray}
f_{1}\left( t_{1}\cdot w_{1}\right) &\leq &\overline{c}_{1}h_{1}\left(
t_{1}\right) \cdot \overline{\varphi }_{1}\left( w_{1}\right) \text{ }%
\forall \text{ }w_{1}\geq 1\text{ and }\forall \text{ }t_{1}\geq M_{1}\cdot
f_{2}^{1/k_{2}}\left( a\right) ,  \label{c21} \\
f_{2}\left( t_{2}\cdot w_{2}\right) &\leq &\overline{c}_{2}h_{2}\left(
t_{2}\right) \cdot \overline{\varphi }_{2}\left( w_{2}\right) \text{ }%
\forall \text{ }w_{2}\geq 1\text{ and }\forall \text{ }t_{2}\geq M_{2}\cdot
f_{1}^{1/k_{1}}\left( b\right) ,  \label{c22}
\end{eqnarray}%
where $M_{1}$ and $M_{2}$ are 
\begin{equation*}
M_{1}=\left\{ 
\begin{array}{lll}
\frac{b}{f_{2}^{1/k_{2}}\left( a\right) } & if & b>f_{2}^{1/k_{2}}\left(
a\right) , \\ 
1 & if & \text{ }b\leq f_{2}^{1/k_{2}}\left( a\right) ,%
\end{array}%
\right. \text{ \ and\ \ }M_{2}=\left\{ 
\begin{array}{lll}
\frac{a}{f_{1}^{1/k_{1}}\left( b\right) } & if & \text{ }a>f_{1}^{1/k_{1}}%
\left( b\right) , \\ 
1 & if & \text{ }a\leq f_{1}^{1/k_{1}}\left( b\right) ,%
\end{array}%
\right.
\end{equation*}%
and $a,b\in \left( 0,\infty \right) $;

(C3)\quad there are some constants $\underline{c}_{1},\underline{c}_{2}\in
\left( 0,\infty \right) $ and the continuous functions $\underline{\varphi }%
_{1}$, $\underline{\varphi }_{2}:\left[ 0,\infty \right) \rightarrow \left[
0,\infty \right) $ such that%
\begin{eqnarray}
f_{1}\left( m_{1}w_{1}\right) &\geq &\underline{c}_{1}\underline{\varphi }%
_{1}\left( w_{1}\right) \text{ }\forall \text{ }w_{1}\geq 1,  \label{c31} \\
f_{2}\left( m_{2}w_{2}\right) &\geq &\underline{c}_{2}\underline{\varphi }%
_{2}\left( w_{2}\right) \text{ }\forall \text{ }w_{2}\geq 1,  \label{c32}
\end{eqnarray}%
where $m_{1}=\min \left\{ b\text{, }f_{2}^{1/k_{2}}\left( a\right) \right\} $
and $m_{2}=\min \left\{ a\text{, }f_{1}^{1/k_{1}}\left( b\right) \right\} $.

\begin{remark}
Examples of functions satisfying (C2) and (C3) are further discussed in the
book of Krasnosel'skii and Rutickii \cite{KR} (see also Gustavsson,
Maligranda and Peetre \cite{GU}).
\end{remark}

Throughout the paper, we use the following notations 
\begin{eqnarray*}
C_{0} &=&(N-1)!/\left[ k_{1}!(N-k_{1})!\right] ,C_{00}=(N-1)!/\left[
k_{2}!(N-k_{2})!\right] \\
G_{1}^{-}\left( \xi \right) &=&\frac{\xi ^{k_{2}-N}}{C_{00}}%
e^{-\int_{0}^{\xi }\frac{1}{C_{00}}t^{k_{2}-1}a_{2}\left( t\right) dt},\text{
}G_{1}^{+}\left( \xi \right) =\xi ^{N-1}e^{\int_{0}^{\xi }\frac{1}{C_{00}}%
t^{k_{2}-1}a_{2}\left( t\right) dt}p_{2}\left( \xi \right) , \\
G_{2}^{-}\left( \xi \right) &=&\frac{\xi ^{k_{1}-N}}{C_{0}}e^{-\int_{0}^{\xi
}\frac{1}{C_{0}}t^{k_{1}-1}a_{1}\left( t\right) dt},\text{ }G_{2}^{+}\left(
\xi \right) =\xi ^{N-1}e^{\int_{0}^{\xi }\frac{1}{C_{0}}t^{k_{1}-1}a_{1}%
\left( t\right) dt}p_{1}\left( \xi \right) , \\
\overline{P}_{1,2}\left( r\right) &=&\int_{0}^{r}[G_{2}^{-}\left( y\right)
\int_{0}^{y}G_{2}^{+}\left( t\right) \overline{\varphi }_{1}\left(
1+\int_{0}^{t}(G_{1}^{-}\left( z\right) \int_{0}^{z}G_{1}^{+}\left( \xi
\right) d\xi )^{\frac{1}{k_{2}}}dz\right) dt]^{1/k_{1}}dy, \\
\underline{P}_{1,2}\left( r\right) &=&\int_{0}^{r}[G_{2}^{-}\left( y\right)
\int_{0}^{y}G_{2}^{+}\left( t\right) \underline{\varphi }_{1}\left(
1+\int_{0}^{t}(G_{1}^{-}\left( z\right) \int_{0}^{z}G_{1}^{+}\left( \xi
\right) d\xi )^{\frac{1}{k_{2}}}dz\right) dt]^{1/k_{1}}dy, \\
\overline{P}_{2,1}\left( r\right) &=&\int_{0}^{r}[G_{1}^{-}\left( y\right)
\int_{0}^{y}G_{1}^{+}\left( t\right) \overline{\varphi }_{2}\left(
1+\int_{0}^{t}(G_{2}^{-}\left( z\right) \int_{0}^{z}G_{2}^{+}\left( \xi
\right) d\xi )^{\frac{1}{k_{1}}}dz\right) dt]^{1/k_{2}}dy, \\
\underline{P}_{2,1}\left( r\right) &=&\int_{0}^{r}[G_{1}^{-}\left( y\right)
\int_{0}^{y}G_{1}^{+}\left( t\right) \underline{\varphi }_{2}\left(
1+\int_{0}^{t}(G_{2}^{-}\left( z\right) \int_{0}^{z}G_{2}^{+}\left( \xi
\right) d\xi )^{\frac{1}{k_{1}}}dz\right) dt]^{1/k_{2}}dy, \\
H_{1,2}\left( r\right) &=&\int_{a}^{r}\frac{1}{h_{1}^{1/k_{1}}\left(
M_{1}f_{2}^{1/k_{2}}\left( t\right) \right) }dt\text{, }H_{2,1}\left(
r\right) =\int_{b}^{r}\frac{1}{h_{2}^{1/k_{2}}\left(
M_{2}f_{1}^{1/k_{1}}\left( t\right) \right) }dt,\text{ } \\
\overline{P}_{1,2}\left( \infty \right) &=&\lim_{r\rightarrow \infty }%
\overline{P}_{1,2}\left( r\right) ,\overline{P}_{2,1}\left( \infty \right)
=\lim_{r\rightarrow \infty }\overline{P}_{2,1}\left( r\right) ,\underline{P}%
_{1,2}\left( \infty \right) =\lim_{r\rightarrow \infty }\underline{P}%
_{1,2}\left( r\right) ,\text{ }\underline{P}_{2,1}\left( \infty \right)
=\lim_{r\rightarrow \infty }\underline{P}_{2,1}\left( r\right) \text{, } \\
H_{1,2}\left( \infty \right) &=&\lim_{s\rightarrow \infty }H_{1,2}\left(
s\right) \text{ and }H_{2,1}\left( \infty \right) =\lim_{s\rightarrow \infty
}H_{2,1}\left( s\right) .
\end{eqnarray*}

Our main results are summarized by the following theorems.

\begin{theorem}
\label{th1}A\textit{ssume that }$H_{1,2}\left( \infty \right) =H_{2,1}\left(
\infty \right) =\infty $ and \textrm{(P1),} hold\textit{. Furthermore, if }$%
f_{1}$ and $f_{2}$ satisfy the hypotheses \textrm{(C1)} and \textrm{(C2)}
then the system (\ref{11}) has at least one positive radial solution 
\begin{equation*}
\left( u_{1},u_{2}\right) \in C^{2}\left( \left[ 0,\infty \right) \right)
\times C^{2}\left( \left[ 0,\infty \right) \right)
\end{equation*}%
with central value in $\left( a,b\right) $. Moreover, the following hold:

1.)\quad If $\overline{P}_{1,2}\left( \infty \right) <\infty $ and $%
\overline{P}_{2,1}\left( \infty \right) <\infty $ then $\lim_{r\rightarrow
\infty }u_{1}\left( r\right) <\infty $ and $\lim_{r\rightarrow \infty
}u_{2}\left( r\right) <\infty .$

2.)\quad If in addition, $f_{1}$ and $f_{2}$ satisfy the hypothesis \textrm{%
(C3)}, $\underline{P}_{1,2}\left( \infty \right) =\infty $ and $\underline{P}%
_{2,1}\left( \infty \right) =\infty $ then%
\begin{equation*}
\lim_{r\rightarrow \infty }u_{1}\left( r\right) =\infty \text{ and }%
\lim_{r\rightarrow \infty }u_{2}\left( r\right) =\infty .
\end{equation*}

3.)\quad If in addition, $f_{2}$ satisfy the condition (\ref{c32}), $%
\overline{P}_{1,2}\left( \infty \right) <\infty $ and $\underline{P}%
_{2,1}\left( \infty \right) =\infty $ then $\lim_{r\rightarrow \infty
}u_{1}\left( r\right) <\infty $ and $\lim_{r\rightarrow \infty }u_{2}\left(
r\right) =\infty .$

4.)\quad If in addition, $f_{1}$ satisfy the condition (\ref{c31}), $%
\underline{P}_{1,2}\left( \infty \right) =\infty $ and $\overline{P}%
_{2,1}\left( \infty \right) <\infty $ then $\lim_{r\rightarrow \infty
}u_{1}\left( r\right) =\infty $ and $\lim_{r\rightarrow \infty }u_{2}\left(
r\right) <\infty .$
\end{theorem}

\begin{theorem}
\label{th2}Assume that the hypothesis \textrm{(P1) }holds. \ Then, the
following hold:

i.)\quad If \textrm{(C1), (C2)}, \textrm{(C3)}, $\overline{P}_{1,2}\left(
\infty \right) <H_{1,2}\left( \infty \right) <\infty $ and $\overline{P}%
_{2,1}\left( \infty \right) <H_{2,1}\left( \infty \right) <\infty $ are
satisfied, then system (\ref{11}) has one positive bounded radial solution 
\begin{equation*}
\left( u_{1},u_{2}\right) \in C^{2}\left( \left[ 0,\infty \right) \right)
\times C^{2}\left( \left[ 0,\infty \right) \right) ,
\end{equation*}%
with central value in $\left( a,b\right) ,$ such that%
\begin{equation*}
\left\{ 
\begin{array}{c}
a+\underline{c}_{1}^{1/k_{1}}\underline{P}_{1,2}\left( r\right) \leq
u_{1}\left( r\right) \leq H_{1,2}^{-1}\left( \overline{c}_{1}^{1/k_{1}}%
\overline{P}_{1,2}\left( r\right) \right) , \\ 
b+\underline{c}_{2}^{1/k_{2}}\underline{P}_{2,1}\left( r\right) \leq
u_{2}\left( r\right) \leq H_{2,1}^{-1}\left( \overline{c}_{2}^{1/k_{2}}%
\overline{P}_{2,1}\left( r\right) \right) .%
\end{array}%
\right. \text{ }
\end{equation*}

ii.)\quad If \textrm{(C1), (C2)}, (\ref{c31}), $H_{1,2}\left( \infty \right)
=\infty $, $\underline{P}_{1,2}\left( \infty \right) =\infty $ and $%
\underline{P}_{2,1}\left( \infty \right) <H_{2,1}\left( \infty \right)
<\infty $ are satisfied, then system (\ref{11}) has one positive radial
solution 
\begin{equation*}
\left( u_{1},u_{2}\right) \in C^{2}\left( \left[ 0,\infty \right) \right)
\times C^{2}\left( \left[ 0,\infty \right) \right) ,
\end{equation*}%
with central value in $\left( a,b\right) $, such that $\lim_{r\rightarrow
\infty }u_{1}\left( r\right) =\infty $ and $\lim_{r\rightarrow \infty
}u_{2}\left( r\right) <\infty $.

iii.)\quad If \textrm{(C1), (C2)}, (\ref{c32}), $\underline{P}_{2,1}\left(
\infty \right) =\infty $, $H_{2,1}\left( \infty \right) =\infty $ and $%
\overline{P}_{1,2}\left( \infty \right) <H_{1,2}\left( \infty \right)
<\infty $ are satisfied, then system (\ref{11}) has one positive radial
solution 
\begin{equation*}
\left( u_{1},u_{2}\right) \in C^{2}\left( \left[ 0,\infty \right) \right)
\times C^{2}\left( \left[ 0,\infty \right) \right) ,
\end{equation*}%
with central value in $\left( a,b\right) $, such that $\lim_{r\rightarrow
\infty }u_{1}\left( r\right) <\infty $ and $\lim_{r\rightarrow \infty
}u_{2}\left( r\right) =\infty $.
\end{theorem}

\begin{remark}
Let%
\begin{equation*}
M_{1}^{+}=\underset{t\in \left[ 0,\infty \right) }{\sup }%
\int_{0}^{t}(G_{1}^{-}\left( z\right) \int_{0}^{z}G_{1}^{+}\left( s\right)
ds)^{\frac{1}{k_{2}}}dz\text{ and }M_{2}^{+}=\underset{t\in \left[ 0,\infty
\right) }{\sup }\int_{0}^{t}(G_{2}^{-}\left( z\right)
\int_{0}^{z}G_{2}^{+}\left( s\right) ds)^{\frac{1}{k_{1}}}dz.
\end{equation*}%
The following situations improve our theorems:

a)\quad If $M_{1}^{+}\in \left( 0,\infty \right) $ then the condition (\ref%
{c21}) is not necessary but $H_{1,2}\left( r\right) $ must be replaced by 
\begin{equation*}
H_{1,2}\left( r\right) =\int_{a}^{r}\frac{1}{f_{1}^{1/k_{1}}\left(
M_{1}\left( 1+M_{1}^{+}\right) f_{2}^{1/k_{2}}\left( t\right) \right) }dt%
\text{.}
\end{equation*}

b)\quad If $M_{2}^{+}\in \left( 0,\infty \right) $ then the condition (\ref%
{c22}) is not necessary but $H_{2,1}\left( r\right) $ must be replaced by 
\begin{equation*}
H_{2,1}\left( r\right) =\int_{b}^{r}\frac{1}{f_{2}^{1/k_{1}}\left(
M_{2}(1+M_{2}^{+})f_{1}^{1/k_{1}}\left( t\right) \right) }dt\text{.}
\end{equation*}

c)\quad If $M_{1}^{+}\in \left( 0,\infty \right) $ and $M_{2}^{+}\in \left(
0,\infty \right) $ then the conditions (\ref{c21}) \ and (\ref{c22}) are not
necessary but $H_{1,2}\left( r\right) $ and $H_{2,1}\left( r\right) $ must
be replaced by 
\begin{equation*}
H_{1,2}\left( r\right) =\int_{a}^{r}\frac{1}{f_{1}^{1/k_{1}}\left(
M_{1}\left( 1+M_{1}^{+}\right) f_{2}^{1/k_{2}}\left( t\right) \right) }dt%
\text{ and }H_{2,1}\left( r\right) =\int_{b}^{r}\frac{1}{f_{2}^{1/k_{1}}%
\left( M_{2}(1+M_{2}^{+})f_{1}^{1/k_{1}}\left( t\right) \right) }dt\text{.}
\end{equation*}

d)\quad If $m_{1}\geq 1$ then $\underline{c}_{1}=1$ and $\underline{\varphi }%
_{1}=f_{1}$.

e)\quad If $m_{2}\geq 1$ then $\underline{c}_{2}=1$ and $\underline{\varphi }%
_{2}=f_{2}$.

f)\quad If $m_{1}\geq 1$ and $m_{2}\geq 1$ then $\underline{c}_{1}=%
\underline{c}_{2}=1$, $\underline{\varphi }_{1}=f_{1}$ and $\underline{%
\varphi }_{2}=f_{2}$.
\end{remark}

\begin{remark}
We note that we can have the following situation 
\begin{equation}
H_{1,2}\left( \infty \right) =\infty \text{ if and only if }H_{2,1}\left(
\infty \right) =\infty ,  \label{in}
\end{equation}%
as the nonlinearities satisfying (\ref{vic}) show. On the other hand we can
construct an example such that (\ref{in}) is not true.
\end{remark}

\begin{remark}
Using the reference \cite{CD} and working as above we can obtain similar
results, as in Theorem \ref{th1} and Theorem \ref{th2}, if the
nonlinearities $f_{1}$, $f_{2}$, $h_{1}$ and $h_{2}$ are assumed to satisfy
the conditions of the Keller-Osserman type%
\begin{equation}
\int_{a}^{\infty }\frac{1}{\sqrt[k_{1}+1]{\left( k_{1}+1\right)
\int_{0}^{s}h_{1}^{1/k_{1}}(c_{1}f_{2}^{1/k_{2}}\left( s\right) )ds}}ds\leq
\infty \text{, }\int_{b}^{\infty }\frac{1}{\sqrt[k_{2}+1]{\left(
k_{2}+1\right) \int_{0}^{s}h_{2}^{1/k_{2}}(c_{2}f_{1}^{1/k_{1}}\left(
s\right) )ds}}ds\leq \infty ,  \label{ko}
\end{equation}%
but the results will not be so strong as here. There are, some differences!
\end{remark}

\section{Proofs of the main results}

In this section we give the proofs of Theorems \ref{th1}-\ref{th2}. \ The
first important tool in our proof is a variant of the Arzel\`{a}--Ascoli
Theorem.

\subsection{The Arzel\`{a}--Ascoli Theorem}

Let $r_{1},r_{2}\in \mathbb{R}$ with $r_{1}\leq r_{2}$ and $\left( K=\left[
r_{1},r_{2}\right] ,d_{K}\left( x,y\right) \right) $ be a compact metric
space, with the metric $d_{K}\left( x,y\right) =\left\vert x-y\right\vert $,
and let 
\begin{equation*}
C\left( \left[ r_{1},r_{2}\right] \right) =\left\{ \dot{g}:\left[ r_{1},r_{2}%
\right] \rightarrow \mathbb{R}\left\vert g\text{ is continuous on }\left[
r_{1},r_{2}\right] \right. \right\}
\end{equation*}%
denote the space of real valued continuous functions on $\left[ r_{1},r_{2}%
\right] $ and for any $g\in C\left( \left[ r_{1},r_{2}\right] \right) $, let 
\begin{equation*}
\left\Vert g\right\Vert _{\infty }=\underset{x\in \left[ r_{1},r_{2}\right] }%
{\max }\left\vert g\left( x\right) \right\vert
\end{equation*}%
be the maximum norm on $C\left( \left[ r_{1},r_{2}\right] \right) $.

\begin{remark}
If $d\left( g^{1}\left( x\right) ,g^{2}\left( x\right) \right) =\left\Vert
g^{1}\left( x\right) -g^{2}\left( x\right) \right\Vert _{\infty }$ then $%
\left( C\left( \left[ r_{1},r_{2}\right] \right) ,d\right) $ is a complete
metric space.
\end{remark}

\begin{definition}
We say that the sequence $\left\{ g_{n}\right\} _{n\in \mathbb{N}}$ from $%
C\left( \left[ r_{1},r_{2}\right] \right) $ is bounded if there exists a
positive constant $C<\infty $ such that $\left\Vert g_{n}\left( x\right)
\right\Vert _{\infty }\leq C$ for each $x\in \left[ r_{1},r_{2}\right] $ (.
Equivalently: $\left\vert g_{n}\left( x\right) \right\vert \leq C$ for each $%
x\in \left[ r_{1},r_{2}\right] $ and $n\in \mathbb{N}^{\ast }$).
\end{definition}

\begin{definition}
We say that the sequence $\left\{ g_{n}\right\} _{n\in \mathbb{N}}$ from $%
C\left( \left[ r_{1},r_{2}\right] \right) $ is equicontinuous if for any
given $\varepsilon >0$, there exists a number $\delta >0$ (which depends
only on $\varepsilon $) such that 
\begin{equation*}
\left\vert g_{n}\left( x\right) -g_{n}\left( y\right) \right\vert
<\varepsilon \text{ for all }n\in \mathbb{N}
\end{equation*}%
whenever $d_{K}\left( x,y\right) <\delta $ for every $x,y\in \left[
r_{1},r_{2}\right] $.
\end{definition}

\begin{definition}
Let $\left\{ g_{n}\right\} _{n\in \mathbb{N}}$ be a family of functions
defined on $\left[ r_{1},r_{2}\right] $. The sequence $\left\{ g_{n}\right\}
_{n\in \mathbb{N}}$ converges uniformly to $g\left( x\right) $ if for every $%
\varepsilon >0$ there is an $N$ (which depends only on $\varepsilon $) such
that%
\begin{equation*}
\left\vert g_{n}\left( x\right) -g\left( x\right) \right\vert <\varepsilon 
\text{ for all }n>N\text{ and }x\in \left[ r_{1},r_{2}\right] .
\end{equation*}
\end{definition}

\begin{theorem}[Arzel\`{a}--Ascoli theorem]
\label{arzela}If a sequence $\left\{ g_{n}\right\} _{n\in \mathbb{N}}$ in $%
C\left( \left[ r_{1},r_{2}\right] \right) $ is bounded and equicontinuous
then it has a subsequence $\left\{ g_{n_{k}}\right\} _{k\in \mathbb{N}}$
which converges uniformly to to $g\left( x\right) $ on $C\left( \left[
r_{1},r_{2}\right] \right) $.
\end{theorem}

\subparagraph{\textbf{Proof of the Theorems \protect\ref{th1} and \protect
\ref{th2}. }}

Recall that 
\begin{eqnarray*}
\lambda \left( D^{2}u\left( x\right) \right) &=&\left\{ 
\begin{array}{l}
\left( \xi ^{\prime \prime }\left( r\right) ,\frac{\xi ^{\prime }\left(
r\right) }{r},...,\frac{\xi ^{\prime }\left( r\right) }{r}\right) \text{ for 
}r\in \left( 0,R\right) , \\ 
\left( \xi ^{\prime \prime }\left( 0\right) ,\xi ^{\prime \prime }\left(
0\right) ,...,\xi ^{\prime \prime }\left( 0\right) \right) \text{ for }r=0%
\end{array}%
\right. \\
S_{k}\left( \lambda \left( D^{2}u\left( x\right) \right) \right) &=&\left\{ 
\begin{array}{l}
C_{N-1}^{k-1}\xi ^{\prime \prime }(r)\left( \frac{\xi ^{\prime }(r)}{r}%
\right) ^{k-1}+C_{N-1}^{k-1}\frac{N-k}{k}\left( \frac{\xi ^{\prime }(r)}{r}%
\right) ^{k}\text{ for }r\in \left( 0,R\right) , \\ 
C_{N}^{k}\left( \xi ^{\prime \prime }\left( 0\right) \right) ^{k}\text{ for }%
r=0,%
\end{array}%
\right.
\end{eqnarray*}%
where for $r=\left\vert x\right\vert <R$, $u\left( x\right) =\xi \left(
r\right) \in C^{2}\left[ 0,R\right) $ is radially symmetric with $\xi
^{\prime }\left( 0\right) =0$ and $C_{N-1}^{k-1}=(N-1)!/\left[ (k-1)!(N-k)!%
\right] $ (see, for example, \cite{BAO} or \cite{SA}).

We start by showing that the system (\ref{11}) has positive radial
solutions. Therefore it can be reduced to 
\begin{equation}
\left\{ 
\begin{array}{l}
C_{N-1}^{k_{1}-1}\left[ \frac{r^{N-k_{1}}}{k_{1}}e^{\int_{0}^{r}\frac{1}{%
C_{0}}t^{k_{1}-1}a_{1}\left( t\right) dt}\left( u_{1}^{^{\prime }}\left(
r\right) \right) ^{k_{1}}\right] ^{\prime }=r^{N-1}e^{\int_{0}^{r}\frac{1}{%
C_{0}}t^{k_{1}-1}a_{1}\left( t\right) dt}p_{1}\left( r\right) f_{1}\left(
u_{2}\left( r\right) \right) \text{ for }r>0, \\ 
C_{N-1}^{k_{2}-1}\left[ \frac{r^{N-k_{2}}}{k_{2}}e^{\int_{0}^{r}\frac{1}{%
C_{00}}t^{k_{2}-1}a_{2}\left( t\right) dt}\left( u_{2}^{^{\prime }}\left(
r\right) \right) ^{k_{2}}\right] ^{\prime }=r^{N-1}e^{\int_{0}^{r}\frac{1}{%
C_{00}}t^{k_{2}-1}a_{2}\left( t\right) dt}p_{2}\left( r\right) f_{2}\left(
u_{1}\left( r\right) \right) \text{ for }r>0, \\ 
u_{1}^{\prime }\left( r\right) \geq 0\text{ and }u_{2}^{\prime }\left(
r\right) \geq 0\text{ for }r\in \left[ 0,\infty \right) , \\ 
u_{1}\left( 0\right) =a\text{ and }u_{2}\left( 0\right) =b.%
\end{array}%
\right.  \label{ss}
\end{equation}
The solutions of (\ref{ss}) can be obtained by using successive
approximation. We define the sequences $\left\{ u_{1}^{m}\right\} ^{m\geq 1}$
and $\left\{ u_{2}^{m}\right\} ^{m\geq 1}$ on $\left[ 0,\infty \right) $ in
the following way: 
\begin{equation}
\begin{array}{l}
u_{1}^{0}=a,u_{2}^{0}=b\text{ for }r\geq 0, \\ 
u_{1}^{m}\left( r\right) =a+\int_{0}^{r}\left[ G_{2}^{-}\left( t\right)
\int_{0}^{t}G_{2}^{+}\left( s\right) f_{1}\left( u_{2}^{m-1}\left( s\right)
\right) ds\right] ^{1/k_{1}}dt, \\ 
u_{2}^{m}\left( r\right) =b+\int_{0}^{r}\left[ G_{1}^{-}\left( t\right)
\int_{0}^{t}G_{1}^{+}\left( s\right) f_{2}\left( u_{1}^{m}\left( s\right)
\right) ds\right] ^{1/k_{2}}dt.%
\end{array}
\label{recs}
\end{equation}%
We can see that $\left\{ u_{1}^{m}\right\} ^{m\geq 1}$ and $\left\{
u_{2}^{m}\right\} ^{m\geq 1}$ are non-decreasing on $\left[ 0,\infty \right) 
$. To do this, let us consider 
\begin{eqnarray*}
u_{1}^{1}\left( r\right) &=&a+\int_{0}^{r}\left[ G_{2}^{-}\left( t\right)
\int_{0}^{t}G_{2}^{+}\left( s\right) f_{1}\left( u_{2}^{0}\left( s\right)
\right) ds\right] ^{1/k_{1}}dt \\
&=&a+\int_{0}^{r}\left[ G_{2}^{-}\left( t\right) \int_{0}^{t}G_{2}^{+}\left(
s\right) f_{1}\left( b\right) ds\right] ^{1/k_{1}}dt \\
&\leq &a+\int_{0}^{r}\left[ G_{2}^{-}\left( t\right)
\int_{0}^{t}G_{2}^{+}\left( s\right) f_{1}\left( u_{2}^{1}\left( s\right)
\right) ds\right] ^{1/k_{1}}dt=u_{1}^{2}\left( r\right) .
\end{eqnarray*}%
This implies that 
\begin{equation*}
u_{1}^{1}\left( r\right) \leq u_{1}^{2}\left( r\right) \text{ which further
produces }u_{2}^{1}\left( r\right) \leq u_{2}^{2}\left( r\right) .
\end{equation*}%
A mathematical induction argument applied to (\ref{recs}) show that for any $%
r\geq 0$ 
\begin{equation*}
u_{1}^{m}\left( r\right) \leq u_{1}^{m+1}\left( r\right) \text{ and }%
u_{2}^{m}\left( r\right) \leq u_{2}^{m+1}\left( r\right) \text{ for any }%
m\in \mathbb{N}\text{,}
\end{equation*}%
i.e.,$\left\{ u_{1}^{m}\right\} ^{m\geq 1}$ and $\left\{ u_{2}^{m}\right\}
^{m\geq 1}$ are non-decreasing on $\left[ 0,\infty \right) $. We now prove
that the non-decreasing sequences $\left\{ u_{1}^{m}\right\} ^{m\geq 1}$ and 
$\left\{ u_{2}^{m}\right\} ^{m\geq 1}$ are bounded from above on bounded
sets. Indeed, by the monotonicity of $\left\{ u_{1}^{m}\right\} ^{m\geq 1}$
and $\left\{ u_{2}^{m}\right\} ^{m\geq 1}$ we have the inequalities%
\begin{eqnarray}
C_{N-1}^{k_{1}-1}\left\{ \frac{r^{N-k_{1}}}{k_{1}}e^{\int_{0}^{r}\frac{1}{%
C_{0}}t^{k_{1}-1}a_{1}\left( t\right) dt}\left[ \left( u_{1}^{m}\left(
r\right) \right) ^{\prime }\right] ^{k_{1}}\right\} ^{\prime }
&=&G_{2}^{+}\left( r\right) f_{1}\left( u_{2}^{m-1}\left( r\right) \right)
\leq G_{2}^{+}\left( r\right) f_{1}\left( u_{2}^{m}\left( r\right) \right) ,
\label{gen1} \\
C_{N-1}^{k_{2}-1}\left\{ \frac{r^{N-k_{2}}}{k_{2}}e^{\int_{0}^{r}\frac{1}{%
C_{00}}t^{k_{2}-1}a_{2}\left( t\right) dt}\left[ \left( u_{2}^{^{m}}\left(
r\right) \right) ^{\prime }\right] ^{k_{2}}\right\} ^{\prime }
&=&G_{1}^{+}\left( r\right) f_{2}\left( u_{1}^{m}\left( r\right) \right) .
\label{gen2}
\end{eqnarray}%
Integrating (\ref{gen1}) leads to 
\begin{eqnarray}
\left( u_{1}^{m}\left( r\right) \right) ^{\prime } &=&\left[ G_{2}^{-}\left(
r\right) \int_{0}^{r}G_{2}^{+}\left( t\right) f_{1}\left( u_{2}^{m-1}\left(
t\right) \right) dt\right] ^{\frac{1}{k_{1}}}  \notag \\
&\leq &\left[ G_{2}^{-}\left( r\right) \int_{0}^{r}G_{2}^{+}\left( t\right)
f_{1}\left( u_{2}^{m}\left( t\right) \right) dt\right] ^{\frac{1}{k_{1}}} 
\notag \\
&=&[G_{2}^{-}\left( r\right) \int_{0}^{r}G_{2}^{+}\left( t\right)
f_{1}\left( b+\int_{0}^{t}(G_{1}^{-}\left( z\right)
\int_{0}^{z}G_{1}^{+}\left( s\right) f_{2}\left( u_{1}^{m}\left( s\right)
\right) ds)^{\frac{1}{k_{2}}}dz\right) dt]^{\frac{1}{k_{1}}}  \notag \\
&\leq &[G_{2}^{-}\left( r\right) \int_{0}^{r}G_{2}^{+}\left( t\right)
f_{1}\left( b+f_{2}^{\frac{1}{k_{2}}}\left( u_{1}^{m}\left( t\right) \right)
\int_{0}^{t}(G_{1}^{-}\left( z\right) \int_{0}^{z}G_{1}^{+}\left( s\right)
ds)^{\frac{1}{k_{2}}}dz\right) dt]^{\frac{1}{k_{1}}}  \notag \\
&\leq &\{G_{2}^{-}\left( r\right) \int_{0}^{r}G_{2}^{+}\left( t\right)
f_{1}\left( f_{2}^{1/k_{2}}\left( u_{1}^{m}\left( t\right) \right) \cdot
\lbrack \frac{b}{f_{2}^{1/k_{2}}\left( u_{1}^{m}\left( t\right) \right) }%
+\int_{0}^{t}(G_{1}^{-}\left( z\right) \int_{0}^{z}G_{1}^{+}\left( s\right)
ds)^{\frac{1}{k_{2}}}dz]\right) dt\}^{\frac{1}{k_{1}}}  \label{exin} \\
&\leq &\{G_{2}^{-}\left( r\right) \int_{0}^{r}G_{2}^{+}\left( t\right)
f_{1}\left( f_{2}^{1/k_{2}}\left( u_{1}^{m}\left( t\right) \right) \cdot
\lbrack \frac{b}{f_{2}^{1/k_{2}}\left( a\right) }+\int_{0}^{t}(G_{1}^{-}%
\left( z\right) \int_{0}^{z}G_{1}^{+}\left( s\right) ds)^{\frac{1}{k_{2}}%
}dz]\right) dt\}^{\frac{1}{k_{1}}}  \notag \\
&\leq &\{G_{2}^{-}\left( r\right) \int_{0}^{r}G_{2}^{+}\left( t\right)
f_{1}\left( M_{1}f_{2}^{1/k_{2}}\left( u_{1}^{m}\left( t\right) \right)
[1+\int_{0}^{t}(G_{1}^{-}\left( z\right) \int_{0}^{z}G_{1}^{+}\left(
s\right) ds)^{\frac{1}{k_{2}}}dz]\right) dt\}^{\frac{1}{k_{1}}}  \notag \\
&\leq &[G_{2}^{-}\left( r\right) \overline{c}_{1}h_{1}\left(
M_{1}f_{2}^{1/k_{2}}\left( u_{1}^{m}\left( r\right) \right) \right)
\int_{0}^{r}G_{2}^{+}\left( t\right) \overline{\varphi }_{1}\left(
1+\int_{0}^{t}(G_{1}^{-}\left( z\right) \int_{0}^{z}G_{1}^{+}\left( s\right)
ds)^{\frac{1}{k_{2}}}dz\right) dt]^{\frac{1}{k_{1}}}  \notag \\
&=&h_{1}^{1/k_{1}}\left( M_{1}f_{2}^{1/k_{2}}\left( u_{1}^{m}\left( r\right)
\right) \right) \overline{c}_{1}^{1/k_{1}}[G_{2}^{-}\left( r\right)
\int_{0}^{r}G_{2}^{+}\left( t\right) \overline{\varphi }_{1}\left(
1+\int_{0}^{t}(G_{1}^{-}\left( z\right) \int_{0}^{z}G_{1}^{+}\left( s\right)
ds)^{\frac{1}{k_{2}}}dz\right) dt]^{\frac{1}{k_{1}}}.  \notag
\end{eqnarray}%
Arguing as above, but now with the second equation (\ref{gen2}), we obtain%
\begin{eqnarray}
\left( u_{2}^{m}\left( r\right) \right) ^{\prime } &=&\left[ G_{1}^{-}\left(
r\right) \int_{0}^{r}G_{1}^{+}\left( z\right) f_{2}\left( u_{1}^{m}\left(
z\right) \right) ds\right] ^{1/k_{2}}  \notag \\
&\leq &h_{2}^{\frac{1}{k_{2}}}\left( M_{2}f_{1}^{\frac{1}{k_{1}}}\left(
u_{2}^{m}\left( r\right) \right) \right) \overline{c}_{2}^{\frac{1}{k_{2}}%
}[G_{1}^{-}\left( r\right) \int_{0}^{r}G_{1}^{+}\left( t\right) \overline{%
\varphi }_{2}\left( 1+\int_{0}^{t}(G_{2}^{-}\left( z\right)
\int_{0}^{z}G_{2}^{+}\left( s\right) ds)^{\frac{1}{k_{1}}}dz\right) dt]^{%
\frac{1}{k_{2}}}.  \label{exin2}
\end{eqnarray}%
Combining the previous relations we obtain 
\begin{eqnarray}
\frac{\left( u_{1}^{m}\left( r\right) \right) ^{\prime }}{h_{1}^{\frac{1}{%
k_{1}}}\left( M_{1}f_{2}^{\frac{1}{k_{2}}}\left( u_{1}^{m}\left( r\right)
\right) \right) } &\leq &\overline{c}_{1}^{\frac{1}{k_{1}}}[G_{2}^{-}\left(
r\right) \int_{0}^{r}G_{2}^{+}\left( t\right) \overline{\varphi }_{1}\left(
1+\int_{0}^{t}(G_{1}^{-}\left( z\right) \int_{0}^{z}G_{1}^{+}\left( s\right)
ds)^{\frac{1}{k_{2}}}dz\right) dt]^{\frac{1}{k_{1}}},  \label{mat1} \\
\frac{\left( u_{2}^{m}\left( r\right) \right) ^{\prime }}{h_{2}^{\frac{1}{%
k_{2}}}\left( M_{2}f_{1}^{\frac{1}{k_{1}}}\left( u_{1}^{m}\left( r\right)
\right) \right) } &\leq &\overline{c}_{2}^{\frac{1}{k_{2}}}[G_{1}^{-}\left(
r\right) \int_{0}^{r}G_{1}^{+}\left( t\right) \overline{\varphi }_{2}\left(
1+\int_{0}^{t}(G_{2}^{-}\left( z\right) \int_{0}^{z}G_{2}^{+}\left( s\right)
ds)^{\frac{1}{k_{1}}}dz\right) dt]^{\frac{1}{k_{2}}}.  \label{mat2}
\end{eqnarray}%
Integrating the inequalities (\ref{mat1}) and (\ref{mat2}), we obtain%
\begin{equation*}
\int_{a}^{u_{1}^{m}\left( r\right) }\frac{1}{h_{1}^{1/k_{1}}\left(
M_{1}f_{2}^{1/k_{2}}\left( t\right) \right) }dt\leq \overline{c}%
_{1}^{1/k_{1}}\overline{P}_{1,2}\left( r\right) \text{ and }%
\int_{b}^{u_{2}^{m}\left( r\right) }\frac{1}{h_{2}^{1/k_{2}}\left(
M_{2}f_{1}^{1/k_{1}}\left( t\right) \right) }dt\leq \overline{c}%
_{2}^{1/k_{2}}\overline{P}_{2,1}\left( r\right) .
\end{equation*}%
We then may write 
\begin{equation}
H_{1,2}\left( u_{1}^{m}\left( r\right) \right) \leq \overline{c}%
_{1}^{1/k_{1}}\overline{P}_{1,2}\left( r\right) \text{ and }H_{2,1}\left(
u_{2}^{m}\left( r\right) \right) \leq \overline{c}_{2}^{1/k_{2}}\overline{P}%
_{2,1}\left( r\right) ,  \label{ints}
\end{equation}%
which plays a basic role in the proof of our main results. Since $H_{\circ
\circ }$ is a bijection with the inverse function $H_{\circ \circ }^{-1}$
strictly increasing on $\left[ 0,\infty \right) $, the inequalities (\ref%
{ints}) can be reformulated as 
\begin{equation}
u_{1}^{m}\left( r\right) \leq H_{1,2}^{-1}\left( \overline{c}_{1}^{1/k_{1}}%
\overline{P}_{1,2}\left( r\right) \right) \text{ and }u_{2}^{m}\left(
r\right) \leq H_{2,1}^{-1}\left( \overline{c}_{2}^{1/k_{2}}\overline{P}%
_{2,1}\left( r\right) \right) .  \label{int}
\end{equation}%
So, we have found upper bounds for 
\begin{equation*}
\left\{ u_{1}^{m}\right\} ^{m\geq 1}\text{ and }\left\{ u_{2}^{m}\right\}
^{m\geq 1}
\end{equation*}%
which are dependent of $r$. We are now ready to give a complete proof of the
Theorems \ref{th1}-\ref{th2}.

\textbf{Proof of Theorem \ref{th1} completed:} We prove that the sequences $%
\left\{ u_{1}^{m}\right\} ^{m\geq 1}$ and $\left\{ u_{2}^{m}\right\} ^{m\geq
1}$ are bounded and equicontinuous on $\left[ 0,c_{0}\right] $ for arbitrary 
$c_{0}>0$. We take 
\begin{equation*}
C_{1}=H_{1,2}^{-1}\left( \overline{c}_{1}^{1/k_{1}}\overline{P}_{1,2}\left(
c_{0}\right) \right) \text{ and }C_{2}=H_{2,1}^{-1}\left( \overline{c}%
_{2}^{1/k_{2}}\overline{P}_{2,1}\left( c_{0}\right) \right)
\end{equation*}%
and since $\left( u_{1}^{m}\left( r\right) \right) ^{^{\prime }}\geq 0$ and $%
\left( u_{2}^{m}\left( r\right) \right) ^{^{\prime }}\geq 0$ it follows that 
\begin{equation*}
u_{1}^{m}\left( r\right) \leq u_{1}^{m}\left( c_{0}\right) \leq C_{1}\text{
and }u_{2}^{m}\left( r\right) \leq u_{2}^{m}\left( c_{0}\right) \leq C_{2}.
\end{equation*}%
We have proved that $\left\{ u_{1}^{m}\left( r\right) \right\} ^{m\geq 1}$
and $\left\{ u_{2}^{m}\left( r\right) \right\} ^{m\geq 1}$ are bounded on $%
\left[ 0,c_{0}\right] $ for arbitrary $c_{0}>0$. Using this fact in (\ref%
{exin}) and (\ref{exin2}) we show that the same is true for $\left(
u_{1}^{m}\left( r\right) \right) ^{\prime }$ and $\left( u_{2}^{m}\left(
r\right) \right) ^{\prime }$. By construction we verify that 
\begin{eqnarray*}
\left( u_{1}^{m}\left( r\right) \right) ^{\prime } &=&\left( \frac{%
r^{k_{1}-N}e^{-\int_{0}^{r}\frac{1}{C_{0}}s^{k_{1}-1}a_{1}\left( s\right) ds}%
}{C_{0}}\int_{0}^{r}t^{N-1}e^{\int_{0}^{t}\frac{1}{C_{0}}s^{k_{1}-1}a_{1}%
\left( s\right) ds}p_{1}\left( t\right) f_{1}\left( u_{2}^{m-1}\left(
t\right) \right) dt\right) ^{1/k_{1}} \\
&\leq &\left( \frac{r^{k_{1}-N}r^{N-1}e^{-\int_{0}^{r}\frac{1}{C_{0}}%
s^{k_{1}-1}a_{1}\left( s\right) ds}}{C_{0}}\int_{0}^{r}p_{1}\left( t\right)
e^{\int_{0}^{t}\frac{1}{C_{0}}s^{k_{1}-1}a_{1}\left( s\right) ds}f_{1}\left(
u_{2}^{m}\left( t\right) \right) dt\right) ^{1/k_{1}} \\
&\leq &\left( \frac{\left\Vert p_{1}\right\Vert _{\infty
}r^{k_{1}-1}e^{-\int_{0}^{r}\frac{1}{C_{0}}s^{k_{1}-1}a_{1}\left( s\right)
ds+\int_{0}^{r}\frac{1}{C_{0}}s^{k_{1}-1}a_{1}\left( s\right) ds}}{C_{0}}%
\int_{0}^{r}f_{1}\left( u_{2}^{m}\left( t\right) \right) dt\right) ^{1/k_{1}}
\\
&\leq &\left( \frac{\left\Vert p_{1}\right\Vert _{\infty }f_{1}\left(
C_{2}\right) r^{k_{1}-1}}{C_{0}}\int_{0}^{r}dt\right) ^{1/k_{1}} \\
&\leq &C_{0}^{-1/k_{1}}\left\Vert p_{1}\right\Vert _{\infty
}^{1/k_{1}}f_{1}^{1/k_{1}}\left( C_{2}\right) c_{0}\text{ on }\left[ 0,c_{0}%
\right] .
\end{eqnarray*}%
Now we turn to $\left( u_{2}^{m}\left( r\right) \right) ^{\prime }$. A
similar argument shows that%
\begin{equation*}
\left( u_{2}^{m}\left( r\right) \right) ^{\prime }\leq
C_{00}^{-1/k_{2}}\left\Vert p_{2}\right\Vert _{\infty
}^{1/k_{2}}f_{2}^{1/k_{2}}\left( C_{1}\right) c_{0}\text{ on }\left[ 0,c_{0}%
\right] .
\end{equation*}%
Finally, it remains to prove that $\left\{ u_{1}^{m}\left( r\right) \right\}
^{m\geq 1}$ and $\left\{ u_{2}^{m}\left( r\right) \right\} ^{m\geq 1}$ are
equicontinuous on $\left[ 0,c_{0}\right] $ for arbitrary $c_{0}>0$. Let $%
\varepsilon _{1}$, $\varepsilon _{2}>0$ be arbitrary. By the mean-value
formula we then deduce that 
\begin{eqnarray*}
\left\vert u_{1}^{m}\left( x\right) -u_{1}^{m}\left( y\right) \right\vert
&=&\left\vert \left( u_{1}^{m}\left( \xi _{1}\right) \right) ^{\prime
}\right\vert \left\vert x-y\right\vert \leq C_{0}^{-1/k_{1}}\left\Vert
p_{1}\right\Vert _{\infty }^{1/k_{1}}f_{1}^{1/k_{1}}\left( C_{2}\right)
c_{0}\left\vert x-y\right\vert , \\
\left\vert u_{2}^{m}\left( x\right) -u_{2}^{m}\left( y\right) \right\vert
&=&\left\vert \left( u_{2}^{m}\left( \xi _{2}\right) \right) ^{\prime
}\right\vert \left\vert x-y\right\vert \leq C_{00}^{-1/k_{2}}\left\Vert
p_{2}\right\Vert _{\infty }^{1/k_{2}}f_{2}^{1/k_{2}}\left( C_{1}\right)
c_{0}\left\vert x-y\right\vert ,
\end{eqnarray*}%
for all $n\in \mathbb{N}$ and all $x,y\in \left[ 0,c_{0}\right] $. So it
suffices to take 
\begin{equation*}
\delta _{1}=\frac{C_{0}^{1/k_{1}}\varepsilon _{1}}{\left\Vert
p_{1}\right\Vert _{\infty }^{1/k_{1}}f_{1}^{1/k_{1}}\left( C_{2}\right) c_{0}%
}\text{ and }\delta _{2}=\frac{C_{00}^{1/k_{2}}\varepsilon _{2}}{\left\Vert
p_{2}\right\Vert _{\infty }^{1/k_{2}}f_{2}^{1/k_{2}}\left( C_{1}\right) c_{0}%
}
\end{equation*}%
to see that $\left\{ u_{1}^{m}\left( r\right) \right\} ^{m\geq 1}$ and $%
\left\{ u_{2}^{m}\left( r\right) \right\} ^{m\geq 1}$ are equicontinuous on $%
\left[ 0,c_{0}\right] $. We now conclude immediately with the help of Arzel%
\`{a}--Ascoli theorem, possibly after passing to a subsequence, that the
sequences$\left\{ u_{j}^{m}\right\} _{j=1,2}^{m\geq 1}$ converges uniformly
to $\left\{ u_{j}\right\} _{j=1,2}$ on $\left[ 0,c_{0}\right] $, in the norm 
$C\left[ 0,c_{0}\right] $. At the end of this process, we conclude, by the
arbitrariness of $c_{0}>0$, that $\left( u_{1},u_{2}\right) $ is a positive
entire solution of system (\ref{11}). The solution constructed in this way
is radially symmetric. Going back to the system (\ref{11}), the radial
solutions of (\ref{ss}) are solutions of the ordinary differential equations
system (\ref{11}). We conclude that radial solutions of (\ref{11}) with $%
u_{1}\left( 0\right) =a,$ $u_{2}\left( 0\right) =b$ satisfy:%
\begin{eqnarray}
u_{1}\left( r\right) &=&a+\int_{0}^{r}\left( G_{2}^{-}\left( y\right)
\int_{0}^{y}G_{2}^{+}\left( t\right) f_{1}\left( u_{2}\left( t\right)
\right) dt\right) ^{1/k_{1}}dy,\text{ }r\geq 0,  \label{eq1} \\
u_{2}\left( r\right) &=&b+\int_{0}^{r}\left( G_{1}^{-}\left( y\right)
\int_{0}^{y}G_{1}^{+}\left( t\right) f_{2}\left( u_{1}\left( t\right)
\right) dt\right) ^{1/k_{2}}dy,\text{ }r\geq 0.  \label{eq2}
\end{eqnarray}%
Next we prove that all four statements hold true.

\textbf{1.):} When $\overline{P}_{1,2}\left( \infty \right) <\infty $ and $%
\overline{P}_{2,1}\left( \infty \right) <\infty ,$ using the same arguments
as (\ref{int}), we find from (\ref{eq1}) and (\ref{eq2}) that%
\begin{equation*}
u_{1}\left( r\right) \leq H_{1,2}^{-1}\left( \overline{c}_{1}^{1/k_{1}}%
\overline{P}_{1,2}\left( \infty \right) \right) <\infty \text{ and }%
u_{2}\left( r\right) \leq H_{2,1}^{-1}\left( \overline{c}_{2}^{1/k_{2}}%
\overline{P}_{2,1}\left( \infty \right) \right) <\infty \text{ for all }%
r\geq 0,
\end{equation*}%
and so $\left( u_{1},u_{2}\right) $ is bounded\textbf{, }which completes the
proof. We next consider:

\textbf{2.):} In the case $\underline{P}_{1,2}\left( \infty \right) =%
\underline{P}_{2,1}\left( \infty \right) =\infty $, we observe that 
\begin{eqnarray}
u_{1}\left( r\right) &=&a+\int_{0}^{r}\left( G_{2}^{-}\left( t\right)
\int_{0}^{t}G_{2}^{+}\left( s\right) f_{1}\left( u_{2}\left( s\right)
\right) ds\right) ^{\frac{1}{k_{1}}}dt  \notag \\
&=&a+\int_{0}^{r}[G_{2}^{-}\left( y\right) \int_{0}^{y}G_{2}^{+}\left(
t\right) f_{1}\left( b+\int_{0}^{t}(G_{1}^{-}\left( z\right)
\int_{0}^{z}G_{1}^{+}\left( s\right) f_{2}\left( u_{1}\left( s\right)
\right) ds)^{\frac{1}{k_{2}}}dz\right) dt]^{\frac{1}{k_{1}}}dy  \notag \\
&\geq &a+\int_{0}^{r}\{G_{2}^{-}\left( y\right) \int_{0}^{y}G_{2}^{+}\left(
t\right) f_{1}\left( [b+f_{2}^{\frac{1}{k_{2}}}\left( a\right)
\int_{0}^{t}(G_{1}^{-}\left( z\right) \int_{0}^{z}G_{1}^{+}\left( s\right)
ds)^{\frac{1}{k_{2}}}dz]\right) dt\}^{\frac{1}{k_{1}}}dy  \label{i1} \\
&\geq &a+\int_{0}^{r}\{G_{2}^{-}\left( y\right) \int_{0}^{y}G_{2}^{+}\left(
t\right) f_{1}\left( m_{1}\cdot \lbrack 1+\int_{0}^{t}(G_{1}^{-}\left(
z\right) \int_{0}^{z}G_{1}^{+}\left( s\right) ds)^{\frac{1}{k_{2}}%
}dz]\right) dt\}^{\frac{1}{k_{1}}}dy  \notag \\
&\geq &a+\int_{0}^{r}[G_{2}^{-}\left( y\right) \int_{0}^{y}G_{2}^{+}\left(
t\right) \underline{c}_{1}\underline{\varphi }_{1}\left(
1+\int_{0}^{t}(G_{1}^{-}\left( z\right) \int_{0}^{z}G_{1}^{+}\left( s\right)
ds)^{\frac{1}{k_{2}}}dz\right) dt]^{\frac{1}{k_{1}}}dy  \notag \\
&=&a+\underline{c}_{1}^{1/k_{1}}\underline{P}_{1,2}\left( r\right) .  \notag
\end{eqnarray}%
The same computations as in (\ref{i1}) yields%
\begin{eqnarray*}
u_{2}\left( r\right) &\geq &b+\underline{c}_{2}^{\frac{1}{k_{2}}%
}\int_{0}^{r}[G_{1}^{-}\left( y\right) \int_{0}^{y}G_{1}^{+}\left( t\right) 
\underline{\varphi }_{2}\left( 1+\int_{0}^{t}(G_{2}^{-}\left( z\right)
\int_{0}^{z}G_{2}^{+}\left( s\right) ds)^{\frac{1}{k_{1}}}dz\right) dt]^{%
\frac{1}{k_{2}}}dy \\
&=&b+\underline{c}_{2}^{1/k_{2}}\underline{P}_{2,1}\left( r\right) ,
\end{eqnarray*}%
and passing to the limit as $r\rightarrow \infty $ in (\ref{i1}) and in the
above inequality we conclude that%
\begin{equation*}
\lim_{r\rightarrow \infty }u_{1}\left( r\right) =\lim_{r\rightarrow \infty
}u_{2}\left( r\right) =\infty ,
\end{equation*}%
which yields the result.

\textbf{3.):} In the spirit of 1.) and 2.) above, we have 
\begin{equation*}
u_{1}\left( r\right) \leq H_{1,2}^{-1}\left( \overline{c}_{1}^{1/k_{1}}%
\overline{P}_{1,2}\left( \infty \right) \right) <\infty \text{ and }%
u_{2}\left( r\right) \geq b+\underline{c}_{2}^{1/k_{2}}\underline{P}%
_{2,1}\left( r\right) .
\end{equation*}%
So, if 
\begin{equation*}
\overline{P}_{1,2}\left( \infty \right) <\infty \text{ and }\underline{P}%
_{2,1}\left( \infty \right) =\infty
\end{equation*}%
we have that 
\begin{equation*}
\lim_{r\rightarrow \infty }u_{1}\left( r\right) <\infty \text{ and }%
\lim_{r\rightarrow \infty }u_{2}\left( r\right) =\infty .
\end{equation*}%
In order, to complete the proof of Theorem \ref{th1} it remains to proceed
to the

\textbf{4.): }In this case, we invoke the proof of 3.). We observe that%
\begin{equation}
u_{1}\left( r\right) \geq a+\underline{c}_{1}^{1/k_{1}}\underline{P}%
_{1,2}\left( r\right) \text{ and }u_{2}\left( r\right) \leq
H_{2,1}^{-1}\left( \overline{c}_{2}^{1/k_{2}}\overline{P}_{2,1}\left(
r\right) \right) .  \label{t2}
\end{equation}%
Our conclusion follows now by letting $r\rightarrow \infty $ in (\ref{t2}).

\textbf{Proof of Theorem \ref{th2} completed: }

\textbf{i.) }Combining (\ref{ints}) and the conditions of the theorem,%
\textbf{\ }we are led to 
\begin{equation*}
\begin{array}{l}
H_{1,2}\left( u_{1}^{m}\left( r\right) \right) \leq \overline{c}%
_{1}^{1/k_{1}}\overline{P}_{1,2}\left( \infty \right) <\overline{c}%
_{1}^{1/k_{1}}H_{1,2}\left( \infty \right) <\infty , \\ 
H_{2,1}\left( u_{2}^{m}\left( r\right) \right) \leq \overline{c}%
_{2}^{1/k_{2}}\overline{P}_{2,1}\left( \infty \right) <\overline{c}%
_{2}^{1/k_{2}}H_{2,1}\left( \infty \right) <\infty .%
\end{array}%
\end{equation*}%
On the other hand, since $H_{\circ }^{-1}$ is strictly increasing on $\left[
0,\infty \right) $, we find that%
\begin{equation*}
u_{1}^{m}\left( r\right) \leq H_{1,2}^{-1}\left( \overline{c}_{1}^{1/k_{1}}%
\overline{P}_{1,2}\left( \infty \right) \right) <\infty \text{ and }%
u_{2}^{m}\left( r\right) \leq H_{2,1}^{-1}\left( \overline{c}_{2}^{1/k_{2}}%
\overline{P}_{2,1}\left( \infty \right) \right) <\infty ,
\end{equation*}%
and then the non-decreasing sequences $\left\{ u_{1}^{m}\left( r\right)
\right\} ^{m\geq 1}$ and $\left\{ u_{2}^{m}\left( r\right) \right\} ^{m\geq
1}$ are bounded above for all $r\geq 0$ and all $m$. Combining these two
facts, we conclude that 
\begin{equation*}
\left( u_{1}^{m}\left( r\right) ,u_{2}^{m}\left( r\right) \right)
\rightarrow \left( u_{1}\left( r\right) ,u_{2}\left( r\right) \right) \text{
as }m\rightarrow \infty
\end{equation*}%
and the limit functions $u_{1}$ and $u_{2}$ are positive entire bounded
radial solutions of system (\ref{11}).\textbf{\ }

\textbf{ii.) and iii.): }For the proof, we follow the same steps and
arguments as in the proof of Theorem \ref{th1}.

\textbf{Acknowledgement. }The author would like to thank to the editors and
reviewers for valuable comments and suggestions which contributed to improve
this article.

\end{document}